\documentclass[a4paper, 11pt]{amsart}

\usepackage[mathscr]{eucal}
\usepackage{amsmath}
\usepackage{amssymb}
\usepackage{amsfonts}
\usepackage{amsthm}

\usepackage{cases}

\usepackage[mathscr]{eucal}
\usepackage{amsmath}
\usepackage{amssymb}
\usepackage{amsfonts}
\usepackage{amsthm}
\usepackage[dvips]{graphics, color}

\theoremstyle{plain}
\newtheorem{theorem}{Theorem}[section]
\newtheorem{corollary}{Corollary}[section]

\theoremstyle{remark}

\newtheorem{remark}{Remark}[section]

\numberwithin{equation}{section}

\def\<{\left<} \def\>{\right>}

\def\proof{\noindent{\it Proof. }}
\def\bea{\begin{eqnarray} }
\def\eea{\end{eqnarray} }
\def\be{\begin{equation} }
\def\ee{\end{equation} }
\def\qed{\ifhmode\unskip\nobreak\fi\ifmmode\ifinner\else\hskip5pt
\fi\fi\hbox{\hskip5 pt \vrule width4 pt height6 pt depth1.5 pt \hskip1pt }}

\begin{document}

\title[]{Special slant surfaces with  non-constant mean curvature  in 2-dimensional complex space forms}
\author[]{Toru Sasahara}
\address{Division of Mathematics, Center for Liberal Arts and Sciences, 
Hachinohe Institute of Technology, 
Hachinohe, Aomori, 031-8501, Japan}
\email{sasahara@hi-tech.ac.jp}

\date{}

\begin{abstract}
{\footnotesize
In the late 1990s, B. Y. Chen introduced the notion of special slant surfaces in  K\"{a}hler surfaces
and classified
non-minimal proper special slant surfaces with constant mean curvature in $2$-dimensional  complex space forms.
In this paper, we completely classify 
proper special slant surfaces with non-constant mean curvature in $2$-dimensional  complex space forms.
 }
\end{abstract}


 \maketitle

 \section{Introduction}
 A surface  in a  K\"{a}hler surface
 is called a slant  surface if  its Wirtinger angle is constant. 
 Holomorphic curves and  Lagrangian surfaces are slant 
 surfaces with constant Wirtinger angle $0$ and $\pi/2$,
 respectively.
  A slant surface is said to be proper if it is neither  a holomorphic curve
   nor  a Lagrangian surface.
There exist infinity many 
proper slant surfaces in $2$-dimensional complex space forms 
(see for example \cite{chenvra, chenvra2, chenvra3}). 

In \cite{chensla}, Chen established a pointwise inequality 
between the squared mean curvature and the Gaussian curvature for
proper slant surfaces in 2-dimensional complex space forms.
 A  proper slant surface satisfies the equality case
   of the inequality identically 
 if and only if the shape operator takes a special form.
  As a generalization of such a surface, 
  Chen   introduced a  new class  of slant surfaces in  K\"{a}hler surfaces,
 named special slant surfaces. This class also 
 includes all slant minimal surfaces and slant
pseudo-umbilical  surfaces.
 Moreover, he \cite{chen3}  classified non-minimal proper special slant surfaces with 
constant mean curvature in
2-dimensional complex space forms.

In this paper, we completely classify 
proper special slant surfaces with non-constant mean curvature in $2$-dimensional  complex space forms.


\section{Preliminaries}
Let $M$ be  
a  surface  of a K\"{a}hler surface $\tilde M^2$ 
with complex structure $J$.
 We denote by $\nabla$ and
$\tilde\nabla$ the Levi-Civita connections on $M$ and 
$\tilde M^2$, and by $h$ and $A_{\xi}$
 the second fundamental form and the shape operator with respect to $\xi$, respectively. 
The
formulas of Gauss and Weingarten are given respectively by
\be
\begin{split}
&\tilde \nabla_XY= \nabla_XY+h(X,Y),\nonumber\\
&\tilde\nabla_X \xi = -A_{\xi}X+D_X\xi
\end{split}
\ee
for vector fields $X$, $Y$ tangent to $M$ and $\xi$ normal to $M$.
The second fundamental form $h$ is related to the shape operator $A$ by
$\<A_{\xi}X, Y\>=\<h(X, Y), \xi\>$.
 The mean curvature vector field $H$ is defined by $H=(1/2){\rm trace}\hskip2pt h$.
The function $|H|$ is called the  {\it mean curvature}.
If it vanishes identically, then $M$ is said to be {\it minimal}.
In particular, if $h$ vanishes identically, then $M$ is said to be {\it totally geodesic}.
A surface $M$ is said to be {\it pseudo-umbilical} if $|H|\ne 0$ and $A_H=\mu I$ for
some function $\mu$ on $M$, where $I$ is the identify transformation.

We denote by $R$ and $\tilde R$ the Riemannian curvature tensors of $M$ and $\tilde M$, respectively,
 and by $R^D$
the curvature tensor of the normal connection $D$. We define
the covariant derivative $\bar\nabla h$ of $h$ by
\be ({\bar\nabla}_{X}h)(Y,Z)= D_X h(Y,Z) - h(\nabla_X
Y,Z) - h(Y,\nabla_X Z).\nonumber
\ee
Then, 
the equations of Gauss, Codazzi and Ricci are  given respectively by
\begin{align}
\<R(X, Y)Z, W\>&=\langle\tilde R(X, Y)Z, W\rangle
+\<h(X, W), h(Y, Z)\>\nonumber\\
&-\<h(X, Z), h(Y, W)\>, \label{ga}\\
(\tilde R(X, Y)Z)^{\perp}&=({\bar\nabla}_{X}h)(Y,Z)-
({\bar\nabla}_{Y}h)(X,Z),\label{co}\\
\langle R^D(X, Y)\xi, \nu\rangle&=\langle \tilde R(X, Y)\xi, \nu\rangle+\<[A_{\xi}, A_{\nu}](X), Y\>
\label{ri}
\end{align}
for vectors $X,Y,Z, W$ tangent  to $M$ and $\xi, \nu$ normal to $M$, where $(\tilde R(X, Y)Z)^{\perp}$
 denotes the normal component of  $\tilde R(X, Y)Z$.

At each point $p\in M$, the angle $\theta$ between $JX$ and $T_pM$, 
where $0\leq\theta\leq\pi/2$, is independent of the choice of non-zero vector $X\in T_pM$.
The angle $\theta$ is called the {\it Wirtinger angle}. 
 A surface is called a {\it slant surface} if its Wirtinger angle
is constant on $M$ (see \cite{chen6}). The Wirtinger angle $\theta$ of a slant surface
is called the {\it slant angle}. A slant surface with slant angle $\theta$ is called $\theta$-{\it slant}. 
Holomorphic curves and 
Lagrangian surfaces are slant surfaces
with slant angle $0$ and $\pi/2$, respectively. 
A slant surface is said to be {\it proper}
if $\theta\ne 0, \pi/2$.

For any vector $X$ tangent to $M$, we put $JX=PX+FX$, where $PX$ and $FX$ denote the tangential and normal 
components of $JX$, respectively. 
If $M$ is a proper slant surface, then for a given unit tangent vector field $e_1$, 
we  choose a local orthonormal frame $\{e_1, e_2, e_3, e_4\}$ of $M$ in  $\tilde M^2$
such that
\be
e_2=(\sec\theta)Pe_1, 
\quad e_3=(\csc\theta)Fe_1, \quad e_4=(\csc\theta)Fe_2.\nonumber
\ee
In the case when $M$ is a Lagrangian surface, for a given local orthonormal frame $\{e_1, e_2\}$ on $M$, 
we choose a local  orthonormal frame
$\{e_1, e_2, e_3, e_4\}$ of $M$ in  $\tilde M^2$ such that 
$e_3=Je_1$ and $e_4=Je_2$.
We call such  orthonormal frames {\it adapted frames}.

\begin{remark}
The equation (\ref{ri})
 of Ricci is a consequence of the equations (\ref{ga}) and (\ref{co})
  of
Gauss and Codazzi for any slant surface with slant angle $\theta\in (0, \pi/2]$
in any K\"{a}hler surface 
(see \cite{chen1}, \cite{chenvra} and \cite{chenvra2}). 
\end{remark}

Let $\tilde M^2(4\epsilon)$ be a complex space form of complex dimension $2$ and constant
holomorphic sectional curvature $4\epsilon$. 
The curvature tensor $\tilde R$ of $\tilde M^2(4\epsilon)$ is given by
\be
\begin{split}
\tilde R(X, Y)Z=&\epsilon\{\<Y, Z\>X-\<X, Z\>Y+\<Z, JY\>JX\\
&-\<Z, JX\>JY+2\<X, JY\>JZ\} \label{space}
\end{split}
\ee
for any tangent vector fields $X$, $Y$ and $Z$ of $\tilde M^2(4\epsilon)$.
It is well-known that any $2$-dimensional complete and simply connected complex space form
is holomorphically isometric  to complex Euclidean plane ${\mathbb C}^2$, 
complex projective plane ${\mathbb C}P^2(4\epsilon)$ or complex hyperbolic plane
 ${\mathbb C}H^2(4\epsilon)$,
according as $\epsilon=0$, $\epsilon>0$ or $\epsilon<0$.

In the next section, we recall the definition of special slant surfaces introduced by Chen \cite{chensla}
and some
results on such slant  surfaces
in $\tilde M^2(4\epsilon)$.

\section{Special slant surfaces}
In \cite{chensla}, 
 Chen proved that the squared mean curvature $|H|^2$ and the Gaussian curvature $K$
 of a proper $\theta$-slant surface  in a complex space form
 $\tilde M^2(4\epsilon)$  satisfy
 the following inequality:
\be \label{ineq}
|H|^2(p)\geq 2K(p)-2(1+3\cos^2\theta)\epsilon
\ee
at each point $p\in M$.
The equality case holds at $p$ 
if and only if, with respect to 
some suitable adapted frame $\{e_1, e_2, e_3, e_4\}$ at $p$,
the shape operator  of $M$ at $p$ takes the following form:
\be 
A_{e_3}= \left(
    \begin{array}{cc}
    3\lambda & 0 \\
     0 & \lambda
    \end{array}
  \right),
  \quad
  A_{e_4}= \left(
    \begin{array}{cc}
    0 & \lambda \\
    \lambda & 0
    \end{array}
  \right).\nonumber
\ee
Moreover, he proved that a non-totally geodesic proper slant 
surface attaining
 equality of the equality (\ref{ineq}) identically exists only in ${\mathbb C}H^2(4\epsilon)$ 
 and satisfies $K\equiv 2\epsilon/3$ and $\theta=\cos^{-1}(1/3)$.
In \cite{chenta}, the explicit representation of such a surface was given.
This surface is named Chen's surface in \cite{ken1}.
A proper slant surface in 
 $\tilde M^2(4\epsilon)$ has non-zero parallel mean curvature
and constant Gaussian curvature if and only if it is locally congruent to Chen's surface (see
\cite{hira, ken1, ken2}).

\begin{remark}
Every Lagrangian surface in  $\tilde M^2(4\epsilon)$
satisfies  the inequality (\ref{ineq}) with $\theta=\pi/2$  
(see \cite{chen4}).
 Chen and Vrancken \cite{chenvran} 
 obtained the explicit representation of Lagrangian surfaces satisfying 
 equality of the equality identically. 
\end{remark}

A slant surface with slant angle $\theta\in (0, \pi/2]$ in a
 K\"{a}hler surface is called {\it special slant} if,
with respect to some suitable adapted frame $\{e_1, e_2, e_3, e_4\}$,
the shape operators of $M$ take the following special form:
\be \label{shape}
A_{e_3}= \left(
    \begin{array}{cc}
    c\lambda & 0 \\
     0 & \lambda
    \end{array}
  \right),
  \quad
  A_{e_4}= \left(
    \begin{array}{cc}
    0 & \lambda \\
    \lambda & 0
    \end{array}
  \right)
\ee
for some constant $c$ and some function $\lambda$ (see \cite{chen3, chensla}).

Every slant minimal surface with slant angle  $\theta\in (0, \pi/2]$ in a  K\"{a}hler surface
is a special slant surface satisfying (\ref{shape}) with $c=-1$ 
(see \cite[Lemma 3.1]{chen2} and \cite[Proposition 3]{chensla}).

Special slant surfaces satisfying (\ref{shape}) with $c=1$ are pseudo-umbilical, 
simply called {\it slumbilical surfaces} (see \cite{chenslu}).
Slumbilical surfaces in 
$\tilde M^2(4\epsilon)$ have been classified in \cite{chenslu}.
In particular, we have 
\begin{theorem}[\cite{chenslu}]\label{T1}
Let $M$ be a proper $\theta$-slant surface  in a complex space form
$\tilde M^2(4\epsilon)$  with non-vanishing mean curvature. If $M$ is a slumbilical surface, then  $\epsilon=0$ and 
one of the following two cases occurs$:$

\noindent $(1)$ $M$ is an open portion of the Euclidean $2$-plane
equipped with  the flat metric
\be
g=e^{-2y\cot\theta}\bigl\{dx^2+(ax+b)^2dy^2\bigr\}\nonumber
\ee
 for some real numbers $a$, $b$ with $a\ne 0$. Moreover, up to rigid
 motions of $\mathbb{C}^2$, the immersion is given by
\begin{align*}
 z(x,y)=&\dfrac{(ax+b)^{1+ia^{-1}\csc\theta}}{a+i\csc\theta}e^{-y\cot\theta}\biggl(\cos\Bigl(\sqrt{1+a^2}y\Bigr)\\
&+i\dfrac{a\cos\theta}{\sqrt{1+a^2}}\sin\Bigl(\sqrt{1+a^2}y\Bigr), 
\dfrac{a\sin\theta+i}{\sqrt{1+a^2}}\sin\Bigl(\sqrt{1+a^2}y\Bigr)\biggr){\rm ;}
\end{align*}

\noindent $(2)$  $M$ is an open portion of the Euclidean $2$-plane
equipped with  the flat metric
\be
g=e^{-2y\cot\theta}(dx^2+b^2dy^2)\nonumber
\ee
  for some positive number $b$. Moreover, up to rigid
 motions of $\mathbb{C}^2$, the immersion is given by
\be 
z(x, y)=b\sin\theta\exp\{ib^{-1}x\csc\theta-y\cot\theta\}(\cos y, \sin y).\nonumber
\ee
\end{theorem}

\begin{remark}\label{rem1}
 Slumbilical surfaces  described
   in (1) (resp. (2)) of Theorem \ref{T1} satisfy (\ref{shape}) with $c=1$ and 
 $\lambda=e^{y\cot\theta}/(ax+b)$ (resp. $\lambda=e^{y\cot\theta}/b$) with respect to an adapted frame such that $e_1=e^{y\cot\theta}\partial/\partial x$ and 
$e_2=\lambda\partial/\partial y$
(cf. Case ($\beta$-$b$) in the proof of Theorem 4.1 of \cite{chenslu}).
In particular, both families of surfaces have non-constant mean curvature.
\end{remark}

For each $c\in(2, 5)$, there exists a special slant immersion 
\be \phi_c: M\Bigl(\frac{4(c-2)}{3(c-1)}\epsilon\Bigr)\rightarrow{\mathbb C}H^2(4\epsilon)\nonumber
\ee
of a simply-connected surface
with constant Gaussian curvature $\frac{4(c-2)}{3(c-1)}\epsilon$
into ${\mathbb C}H^2(4\epsilon)$ with slant angle
$\theta=\cos^{-1}\Bigl(\frac{1}{3}\sqrt{\frac{(5-c)(c-2)}{c-1}}\Bigr)$, which satisfies (\ref{shape}) with 
$\lambda=\sqrt{\frac{(c-5)\epsilon}{3(c-1)}}$, with respect to 
some adapted frame (see  \cite[Theorem1]{chen3}).
The explicit representation of $\phi_c$
 was obtained in  \cite[Theorem 5.1]{chenta}.  
A  non-totally geodesic proper slant surface satisfies 
 equality of the equality (\ref{ineq}) identically if and only if 
 it is locally congruent to $\phi_3$. 

Let $U$ be a simply-connected open subset of 
${\mathbb R}^2$ and $E=E(x,y)$ a positive function on $U$ satisfying
the following conditions:
\be
\frac{\partial}{\partial x}\biggl(\frac{1}{E}\frac{\partial}
{\partial x}\Bigl(\frac{1}{E}\frac{\partial E}{\partial y}\Bigr)\biggr)=0,\quad
\frac{\partial E}{\partial y}\ne 0.\nonumber
\ee
For a given $\theta\in(0, \pi/2)$ and $\epsilon<0$,
 we put $G=(\sec\theta)E_y/(2\sqrt{-\epsilon}E)$, 
$E_y=\partial E/\partial y$.
Denote by $M(\theta, \epsilon, E)$ the surface $U$ with metric tensor
\be
g=E^2dx^2+G^2dy^2.\nonumber
\ee
Then,  $M(\theta, \epsilon, E)$  has constant Gaussian curvature  $K=4\epsilon\cos^2\theta$.
There exists a  special $\theta$-slant immersion 
 \be 
 \psi_{\theta, \epsilon, E}: M(\theta, \epsilon, E)\rightarrow {\mathbb C}H^2(4\epsilon),\nonumber
 \ee 
 which satisfies (\ref{shape}) with $c=2$ and $\lambda=-\sqrt{-\epsilon}\sin\theta$,
 with respect to the adapted frame with $e_1=E^{-1}\partial/{\partial x}$, 
$e_2=G^{-1}\partial/{\partial y}$ (see \cite[Theorem 2]{chen3}).

The explicit representation of $\psi_{\theta, -1, y}$  was given in Proposition 6.1 
of \cite{chenta}.
Two slant immersions $\psi_{\theta, -1, (x+y)^{-2}}$ and 
$\psi_{\theta,-1, y}$ are not congruent to each other in  ${\mathbb C}H^2(-4)$.
Therefore, up to rigid motions of ${\mathbb C}H^2(-4)$, for each $\theta\in(0, \pi/2)$
there exist more than one special $\theta$-slant  immersions of a surface
with constant Gaussian curvature 
$-4\cos^2\theta$ into ${\mathbb C}H^2(-4)$ satisfying (\ref{shape})
(see Remark 6.1 and Corollary 6.2 in \cite{chenta}).

For special slant surfaces with satisfying (\ref{shape})
with $c=2$ in $\tilde M^2(4\epsilon)$, we have
\begin{theorem}[\cite{chen3}]\label{L3}
Let $M$ be a non-minimal proper $\theta$-slant surface  in a complex space form
$\tilde M^2(4\epsilon)$. If $M$ is  a special slant surface satisfying  $(\ref{shape})$ with $c=2$, then
$\epsilon<0$ and  $M$ is isometric to $M(\theta, \epsilon, E)$. Moreover, up to rigid motions of  
${\mathbb C}H^2(4\epsilon)$, the immersion is given by  
$\psi_{\theta, \epsilon, E}$.
\end{theorem}

Chen \cite{chen3} proved that $\phi_c$ and  $\psi_{\theta, \epsilon, E}$
are the only non-minimal proper special slant 
surfaces in  a complex space form $\tilde M^2(4\epsilon)$ with
constant mean curvature.
The next section solves the classification problem for the case of non-constant mean curvature.


\section{The main theorem}

The main result of this paper is the following classification theorem.
\begin{theorem}\label{main}
Let $M$ be a proper $\theta$-slant surface in a complex space form
$\tilde M^2(4\epsilon)$ 
with $\nabla |H|^2\ne 0$ everywhere.
If $M$ is a special slant surface, then $M$ is a slumbilical surface in $\mathbb{C}^2$, which
 is locally congruent to one of the following$:$

\noindent $(1)$
\begin{align*}
 z(x,y)=&\dfrac{(ax+b)^{1+ia^{-1}\csc\theta}}{a+i\csc\theta}e^{-y\cot\theta}\biggl(\cos\Bigl(\sqrt{1+a^2}y\Bigr)\\
&+i\dfrac{a\cos\theta}{\sqrt{1+a^2}}\sin\Bigl(\sqrt{1+a^2}y\Bigr), 
\dfrac{a\sin\theta+i}{\sqrt{1+a^2}}\sin\Bigl(\sqrt{1+a^2}y\Bigr)\biggr) 
\end{align*}
for some real numbers $a$, $b$ with $a\ne 0${\rm;}

\noindent $(2)$
\be 
z(x, y)=b\sin\theta\exp\{ib^{-1}x\csc\theta-y\cot\theta\}(\cos y, \sin y) \nonumber
\ee
for some positive number $b$.
\end{theorem}

\proof
Let $M$ be a  proper special slant surface in  
$\tilde M^2(4\epsilon)$ 
satisfying (\ref{shape}) with respect to some suitable
adapted frame $\{e_1, e_2, e_3, e_4\}$. Then, we have
\be
h(e_1, e_1)=c\lambda e_3, \ \ h(e_1, e_2)=\lambda e_4, \ \ h(e_2, e_2)=\lambda e_3. \label{second}
\ee
Assume that $\nabla |H|^2\ne 0$ everywhere.
 Then, $c\ne -1$,
$\lambda$ is not constant and nowhere zero.
Moreover, it follows from Theorem \ref{L3} that  $c\ne -2$.

We put $p=\<\nabla_{e_1}e_1, e_2\>$ and $q=\<\nabla_{e_2}e_1, e_2\>$.
Then, by Lemma 4 of \cite{chensla}, which is a consequence of  the equation (\ref{co})
of Codazzi and Lemma 4.1 of \cite[p.29]{chen6}, we have 
\begin{align}
& e_1\lambda=(c-2)\lambda q, \label{e1}\\
& e_2\lambda= \lambda p+\dfrac{6\epsilon}{c+1}\sin\theta\cos\theta,\label{e2}\\
& \lambda p=\dfrac{(c+1)\lambda^2\cot\theta}{2}-\dfrac{3(c-1)}{2(c+1)}\epsilon\sin\theta\cos\theta.\label{e3}
\end{align}
It follows from (\ref{space}), (\ref{second}) and  the equation (\ref{ga}) of Gauss for $\<R(e_1, e_2)e_2, e_1\>$ that
\be
-p^2-q^2-e_1q+e_2p=(c-1)\lambda^2+(1+3\cos^2\theta)\epsilon.\label{e4}
\ee
Differentiating (\ref{e3}) along $e_i$ $(i=1, 2)$ and using (\ref{e1}) and (\ref{e2}), 
we obtain
\begin{align}
& e_1p=(c-2)\bigl\{(c+1)\lambda\cot\theta-p\bigr\}q,\label{e5}\\
& e_2p=-p^2+ (c+1)\lambda p\cot\theta-\dfrac{6\epsilon p}{(c+1)\lambda}\sin\theta\cos\theta+6\epsilon\cos^2\theta.\label{e6}
\end{align}
Substituting (\ref{e6}) into (\ref{e4}) yields
\be 
e_1q=-(c-1)\lambda^2-\epsilon-2p^2-q^2+ (c+1)\lambda p\cot\theta-\dfrac{6\epsilon p}{(c+1)\lambda}\sin\theta\cos\theta+3\epsilon\cos^2\theta.\label{e7}
\ee

Every function $f$ satisfies  the condition 
 $[e_1, e_2]f-(\nabla_{e_1}e_2-\nabla_{e_2}e_1)f=0$,
which can be written as
\be
e_1(e_2f)-e_2(e_1f)+pe_1f+qe_2f=0. \label{integral}
\ee
We apply (\ref{integral}) for $\lambda$. 
Using (\ref{e1}) and (\ref{e2}),  we get
\be
e_1(e_2\lambda)-e_2(e_1\lambda)=(e_1\lambda)p+\lambda e_1p-(c-2)\bigl\{(e_2\lambda)q+\lambda e_2q\bigr\}.\label{integc}
\ee
From (\ref{integral}) with $f=\lambda$ and (\ref{integc}), it follows  that
\be
2pe_1\lambda+\lambda e_1p+(3-c)qe_2\lambda-(c-2)\lambda e_2q=0.\label{integc2}
\ee 
Substituting (\ref{e1}), (\ref{e2}) and (\ref{e5}) into (\ref{integc2}) gives
\be
e_2q=q\left(\dfrac{p}{c-2}+(c+1)\lambda\cot\theta+\dfrac{6(3-c)\epsilon}{(c-2)(c+1)\lambda}\sin\theta\cos\theta\right).\label{e8}
\ee

{\bf Case (I):} $q=0$ on an open subset $\mathcal{U}$ of $M$. 
In this case, substituting (\ref{e3}) into {(\ref{e7}), and 
multiplying the resulting equation by the least common denominator $2(c+1)^2\lambda^2$,
 we obtain the following 
polynomial equation in  $\lambda$: 
\be 2 (c-1) (c+1)^2\lambda^4
-\epsilon (c+1)^2 \{3(c-1)\cos^2\theta-2\}\lambda^2
+9 \epsilon^2 (c-3) (c-1) \sin^2\theta \cos^2\theta=0.\nonumber
\ee
Since $\lambda$ is not constant, all coefficients must vanish. Therefore, taking into
account $c\ne -1$, we must have
 $c=1$. 
Applying Theorem \ref{T1} shows $\epsilon=0$. 
It follows from  $q=0$ and ({\ref{e1}) that $e_1\lambda=0$.
Thus, Remark \ref{rem1} implies that
 $\mathcal{U}$ is locally congruent to (2).

{\bf Case (II):} $q \ne 0$ on an open subset $\mathcal{W}$ of $M$. In this case, 
we compute (\ref{integral}) with $f=q$, that is, 
 \be
 e_1(e_2q)-e_2(e_1q)+pe_1q+qe_2q=0. \label{integq}
 \ee
Differentiating (\ref{e7}) and (\ref{e8}) along $e_2$ and $e_1$ respectively, we obtain
\begin{align}
e_2(e_1q)=&\left\{-2(c-1)\lambda+p(c+1)\cot\theta
+\dfrac{6\epsilon p}{(c+1)\lambda^2}\sin\theta\cos\theta\right\}e_2\lambda \nonumber\\
&-\left\{4p-(c+1)\lambda\cot\theta+\dfrac{6\epsilon}{(c+1)\lambda}\sin\theta\cos\theta\right\}e_2p-2qe_2q,\label{q1}\\
e_1(e_2q)=&\left\{\dfrac{p}{c-2}+(c+1)\lambda\cot\theta+\dfrac{6(3-c)\epsilon}{(c-2)(c+1)\lambda}\sin\theta\cos\theta\right\}e_1q\nonumber\\
&+\dfrac{q}{c-2}e_1p+q\left\{(c+1)\cot\theta-\dfrac{6(3-c)\epsilon}{(c-2)(c+1)\lambda^2}\sin\theta\cos\theta\right\}e_1\lambda.\label{q2}
\end{align}
 Substituting  (\ref{q1}) and (\ref{q2}) into (\ref{integq}) shows
 \be
 \begin{split}
 &\dfrac{q}{c-2}e_1p+\left\{4p-(c+1)\lambda\cot\theta+\dfrac{6\epsilon}{(c+1)\lambda}\sin\theta\cos\theta\right\}e_2p\\
 &+\left\{\dfrac{(c-1)p}{c-2}+(c+1)\lambda\cot\theta+\dfrac{6(3-c)\epsilon}{(c-2)(c+1)\lambda}\sin\theta\cos\theta\right\}e_1q+3qe_2q\\
 &+q\left\{(c+1)\cot\theta-\dfrac{6(3-c)\epsilon}{(c-2)(c+1)\lambda^2}\sin\theta\cos\theta\right\}e_1\lambda\\
 &+ \left\{2(c-1)\lambda-(c+1)p\cot\theta
-\dfrac{6\epsilon p}{(c+1)\lambda^2}\sin\theta\cos\theta\right\}e_2\lambda=0. \label{q3}
\end{split} 
 \ee
 Substituting  (\ref{e1}), (\ref{e2}), (\ref{e5}), (\ref{e6}), (\ref{e7}) and (\ref{e8}) into  (\ref{q3})
and  multiplying  the resulting equation by the least common denominator $(c-2)(c+1)^2\lambda^2$,
we find
\be 
\alpha_1(\lambda)p^3+\alpha_2(\lambda)p^2+\alpha_3(\lambda)p
+\alpha_4(\lambda)=0,\label{int1}
\ee
where $\alpha_i$ $(i=1, 2, 3, 4)$ are 
polynomials in $\lambda$ given by
 \begin{align*}
 \alpha_1=&2  (c+1)^2 (3 c-5)\lambda^2,\\
 \alpha_2=&-(c+1)^3 (3 c-5) \cot\theta \lambda^3+6 \epsilon (c+1) (5 c-7) \sin\theta \cos\theta \lambda,\\
 \alpha_3=&-(c-3) (c-1) (c+1)^2\lambda^4\\
 &+\bigl\{(c+1)^2 [\epsilon(c-1)+2 q^2(c-3)]-3 \epsilon (c+1)^2 (7 c-11) \cos^2\theta\bigr\}\lambda^2\\
 &+36 \epsilon^2 (c-1) \sin^2\theta \cos^2\theta,\\
 \alpha_4=&(c-2) (c-1) (c+1)^3 \cot\theta \lambda^5\\
 &+\bigl\{(c-2) (c+1)^3  [\epsilon-q^2(c+1)]\cot\theta\\
 &+3\epsilon(c-2)(c+1)^3\cos^2\theta\cot\theta-6 \epsilon (c-1) (c+1) (3 c-7) \sin\theta \cos\theta\bigr\}\lambda^3\\
 &-\bigl\{6 \epsilon (c-3) (c+1) [\epsilon+q^2(c-4)]\sin\theta \cos\theta +18 \epsilon^2 (c+1)(c-1) \sin\theta \cos^3\theta\bigr\}\lambda,
 \end{align*}
where all coefficients are
 expressed as \be
\sum_{m, n\geq 1}Q_{m, n}^i(c)\sin^m\theta\cos^n\theta+\sum_{m, n\geq 0}R_{m, n}^i(c)
\cos^m\theta\cot^n\theta \label{coef}
\ee 
 for some polynomials $Q_{m, n}^i(c)$ and $R_{m, n}^i(c)$ in $c$.

 In the rest of this paper, we will express the coefficients of polynomials in $\lambda$   in the 
  form (\ref{coef}). 

Substituting (\ref{e3}) into (\ref{int1}) and multiplying the resulting equation by
 the least common denominator $4(c+1)\lambda$, we derive
 \be
 \beta_1(\lambda)q^2+\beta_2(\lambda)=0,\label{int2}
 \ee
where $\beta_i$ $(i=1, 2)$ are polynomials in  $\lambda$ given by
\begin{align*}
\beta_1=&4 (c-1)^2 (c+1)^4 \cot\theta \lambda^4+36 \epsilon (c-3)^2 (c+1)^2 \sin\theta \cos\theta \lambda^2,\\
\beta_2=&-2 (c-1)^2 (c+1)^4 \cot\theta \lambda^6\\
&-\bigl\{2 \epsilon (c+1)^4 (3 c-5) \cot\theta-3 \epsilon
 (c-1)(c+1)^4 (3 c-5) \cos^2\theta \cot\theta\\
 &+6 \epsilon (c-1) (c+1)^2 (c^2-16 c+31) \sin\theta \cos\theta \bigr\}\lambda^4\\ 
 &+\bigl\{6 \epsilon^2 (c+1)^2 (c^2+2 c-11) \sin\theta \cos\theta\\
 &-18 \epsilon^2 (c-1)^2(c+1)^2(3c-8)\sin\theta \cos^3\theta\bigr\}\lambda^2\\
& +81 \epsilon^3 (c-3)^2 (c-1)^2 \sin^3\theta \cos^3\theta.
\end{align*}
We differentiate (\ref{int2}) along $e_1$ and use
(\ref{e1}) and (\ref{e7}). Then, dividing the  left-hand side of
the resulting equation by the
greatest common factor $4(c+1)\lambda q$, we obtain  
\be 
\kappa_1(\lambda)p^2+
\kappa_2(\lambda)p+\kappa_3(\lambda)=0,\label{int3}
\ee
where $\kappa_i$ $(i=1, 2, 3)$ are polynomials in $\lambda$ given by
\begin{align*}
\kappa_1=&4 (c-1)^2 (c+1)^3 \cot\theta \lambda^3+36 \epsilon (c-3)^2 (c+1) \sin\theta \cos\theta \lambda,\\
\kappa_2=&-2 (c-1)^2 (c+1)^4 \cot^2\theta \lambda^4-6 \epsilon (c+1)^2 (c^2-14 c+25) \cos^2\theta
\lambda^2\\
&+108 \epsilon^2 (c-3)^2 \sin^2\theta \cos^2\theta,\\
\kappa_3=&(c-1)^2 (c+1)^3 (5 c-8) \cot\theta \lambda^5\\
&+\bigl\{2 (c+1)^3  \bigl[\epsilon (4 c^2-13 c+11)-q^2 (c-1)^2 (2 c-5)\bigr]\cot\theta\\
&-3 \epsilon (c-1)(c+1)^3(3c^2-7c+6) \cos^2\theta \cot\theta\\
&+
6 \epsilon (c-1) (c+1) (c^3-15 c^2+45 c-35) \sin\theta \cos\theta\bigr\}\lambda^3\\
&-\bigl\{3 \epsilon (c+1)  \bigl[\epsilon (c^3-6 c^2+21 c-32)+6 q^2 (c-3)^3\bigr]\sin\theta \cos\theta\\
&+9 \epsilon^2 (c+1) (3 c^4-20 c^3+35 c^2+26 c-92) \sin\theta \cos^3\theta\bigr\}\lambda.
\end{align*}
Substituting 
(\ref{e3}) into (\ref{int3}) and multiplying the resulting equation by the least common denominator $(c+1)\lambda$, we get
\be
\tau_1(\lambda)q^2+\tau_2(\lambda)=0, \label{int4}
\ee
where $\tau_i$ $(i=1, 2)$ are polynomials in $\lambda$ given by
\begin{align*}
\tau_1=&2 (c-1)^2 (c+1)^4 (2 c-5) \cot\theta \lambda^4
+18 \epsilon (c-3)^3 (c+1)^2 \sin\theta \cos\theta \lambda^2,\\
\tau_2=&-(c-1)^2 (c+1)^4 (5 c-8) \cot\theta \lambda^6\\
&-\bigl\{2 \epsilon (c+1)^4 (4 c^2-13 c+11) \cot\theta-3 \epsilon (c-1)(c+1)^4(4c^2-11c+9) \cos^2\theta
 \cot\theta\\
 &+6 \epsilon (c-1) (c+1)^2 (c^3-15 c^2+45 c-35) \sin\theta \cos\theta\bigr\}\lambda^4\\
 &+\bigl\{3 \epsilon^2 (c+1)^2 (c^3-6 c^2+21 c-32) \sin\theta \cos\theta\\
 &-9 \epsilon^2 (c+1)^2 (4 c^4-29 c^3+74 c^2-65 c-8) \sin\theta \cos^3\theta\bigr\}\lambda^2\\
& -81 \epsilon^3 (c-3)^3 (c-1) \sin^3\theta \cos^3\theta.
\end{align*}
Eliminating $q^2$ from (\ref{int2}) and (\ref{int4}) yields
\be
\beta_1\tau_2-\beta_2\tau_1=0,\nonumber
\ee
which can be written as the following polynomial equation in $\lambda$:
\be 
12 (c-1)^5 (c+1)^8 \cot^2\theta \lambda^{10}+\epsilon\sum_{m=1}^4\psi_{m}(c)\lambda^{2m}=0,\nonumber
\ee
where each $\psi_{m}(c)$ is a polynomial in $c$.
Since $\lambda$ is not constant, all coefficients must vanish. Hence, taking into account $c\ne -1$, we must have 
$c=1$. Applying Theorem \ref{T1}, we get $\epsilon=0$. 
It follows from  $q\ne 0$, $c\ne 2$ and ({\ref{e1}) that $e_1\lambda\ne 0$.
Therefore, from Remark \ref{rem1} we deduce  that
$\mathcal{W}$ is locally congruent to (1). 
\qed

\medskip
Combining Theorem \ref{main}, Corollary 3.2 of \cite{chenta} and Theorem 5 of \cite{chen3},
we  have
\begin{corollary}\label{cor}
Let $M$ be  a non-totally geodesic proper slant surface in  a complex space form
$\tilde M^2(4\epsilon)$. If $M$ is a  special slant surface
 satisfying $(\ref{shape})$, then either
 
\noindent  $(1)$ $\epsilon=0$ and $c\in\{-1, 1\}$, or

\noindent $(2)$ $\epsilon<0$ and $2\leq c<5$.
\end{corollary}

\section{Final remarks}

\begin{remark}
In contrast to  Corollary \ref{cor},
 for any given constants $\epsilon$ and $c$, there exists a Lagrangian surface
in  
$\tilde M^2(4\epsilon)$  satisfying $(\ref{shape})$ for some non-zero function $\lambda$ with respect to some
adapted frame (see Theorem 6.1 of \cite{chen2}).
\end{remark}

\begin{remark}
Let $M$ be a surface 
in a Riemannian $4$-manifold.
We consider an arbitrary orthonormal basis $\{e_1, e_2\}$ of $T_pM$ and  denote $h_{ij}=h(e_i, e_j)$.  
Then,  for $X=\cos\theta e_1+\sin\theta e_2$ we have 
\be
h(X, X)=H+\cos 2\theta\dfrac{h_{11}-h_{22}}{2}+\sin 2\theta h_{12}.\nonumber
\ee 
This  implies that
when $X$ goes once around the unit tangent circle, $h(X, X)$ goes twice around 
the ellipse centered at $H$, which is called  the {\it curvature ellipse}.
The ellipse can degenerate into a line segment or a point.

For a non-minimal special slant surface satisfying (\ref{shape}) in 
a K\"{a}hler surface, we have
 
\noindent (1) $c=3$ if and only if the curvature ellipse is a circle at any point
(see \cite{cas} for Lagrangian surfaces in complex space forms).

\noindent (2) $c=1$ if and only if the curvature ellipse degenerates into a line 
segment at any point.

\noindent (3)  If the ambient space is a complex space form and 
$c\ne 1$, then the semi-major axis and the semi-minor axis
 of the curvature ellipse are both positive constant.
\end{remark}

\begin{remark}
Theorem 6 in \cite{chensla} states that a proper slant surface in the complex Euclidean
plane $\mathbb{C}^2$ is special slant
if and only if it is a slant minimal surface.
However, it is incorrect. In fact, special slant surfaces in $\mathbb{C}^2$ described 
in Theorem \ref{main} are 
non-minimal surfaces.
In  the proof of Theorem 6 in \cite{chensla}, the case 
$c=1$   is missing.
If $c=1$, then the solution $\lambda$ of  (4.34) in
\cite{chensla} is given by 
 $$\lambda=f(x)e^{\cot\theta\int\phi(y)dy}$$ 
 for some function $f(x)$. 
 If $c\ne 1$, then $\lambda$ is given by (4.35)  in \cite{chensla}, which leads to a contradiction.
  Hence, Theorem 6 in \cite{chensla} should be corrected as follows:
  
  {\it
A  proper slant surface in the complex Euclidean plane
$\mathbb{C}^2$ is special slant
if and only if it is either a slant minimal surface or a slumbilical surface.}
\end{remark}

\begin{remark}
In the proof of Theorem \ref{main},  the condition (\ref{integq})
for $q$
plays  an important role.
On the contrary, 
(\ref{integral}) for $f=p$
does not provide any new information. 
In fact,
 using (\ref{e1}), (\ref{e2}), (\ref{e5}) and (\ref{e6}), we can check
that   (\ref{integral}) for $f=p$ holds identically. 
\end{remark}

\end{document}